\hfuzz 50pt
\font\ninerm=cmr9  \font\eightrm=cmr8  \font\sixrm=cmr6
\font\ninei=cmmi9  \font\eighti=cmmi8  \font\sixi=cmmi6
\font\ninesy=cmsy9 \font\eightsy=cmsy8 \font\sixsy=cmsy6
\font\ninebf=cmbx9 \font\eightbf=cmbx8 \font\sixbf=cmbx6
\font\nineit=cmti9 \font\eightit=cmti8 
\font\ninett=cmtt9 \font\eighttt=cmtt8 
\font\ninesl=cmsl9 \font\eightsl=cmsl8

\font\twelverm=cmr12 at 15pt
\font\twelvei=cmmi12 at 15pt
\font\twelvesy=cmsy10 at 15pt
\font\twelvebf=cmbx12 at 15pt
\font\twelveit=cmti12 at 15pt
\font\twelvett=cmtt12 at 15pt
\font\twelvesl=cmsl12 at 15pt
\font\twelvegoth=eufm10 at 15pt

\font\tengoth=eufm10  \font\ninegoth=eufm9
\font\eightgoth=eufm8 \font\sevengoth=eufm7 
\font\sixgoth=eufm6   \font\fivegoth=eufm5
\newfam\gothfam \def\goth{\fam\gothfam\tengoth} 
\textfont\gothfam=\tengoth
\scriptfont\gothfam=\sevengoth 
\scriptscriptfont\gothfam=\fivegoth

\catcode`@=11
\newskip\ttglue

\def\tenpoint{\def\rm{\fam0\tenrm}
  \textfont0=\tenrm \scriptfont0=\sevenrm
  \scriptscriptfont0\fiverm
  \textfont1=\teni \scriptfont1=\seveni
  \scriptscriptfont1\fivei 
  \textfont2=\tensy \scriptfont2=\sevensy
  \scriptscriptfont2\fivesy 
  \textfont3=\tenex \scriptfont3=\tenex
  \scriptscriptfont3\tenex 
  \textfont\itfam=\tenit\def\it{\fam\itfam\tenit}%
  \textfont\slfam=\tensl\def\sl{\fam\slfam\tensl}%
  \textfont\ttfam=\tentt\def\tt{\fam\ttfam\tentt}%
  \textfont\gothfam=\tengoth\scriptfont\gothfam=\sevengoth 
  \scriptscriptfont\gothfam=\fivegoth
  \def\goth{\fam\gothfam\tengoth}
  \textfont\bffam=\tenbf\scriptfont\bffam=\sevenbf
  \scriptscriptfont\bffam=\fivebf
  \def\bf{\fam\bffam\tenbf}%
  \tt\ttglue=.5em plus.25em minus.15em
  \normalbaselineskip=12pt \setbox\strutbox\hbox{\vrule
  height8.5pt depth3.5pt width0pt}%
  \let\big=\tenbig\normalbaselines\rm}

\def\ninepoint{\def\rm{\fam0\ninerm}
  \textfont0=\ninerm \scriptfont0=\sixrm
  \scriptscriptfont0\fiverm
  \textfont1=\ninei \scriptfont1=\sixi
  \scriptscriptfont1\fivei 
  \textfont2=\ninesy \scriptfont2=\sixsy
  \scriptscriptfont2\fivesy 
  \textfont3=\tenex \scriptfont3=\tenex
  \scriptscriptfont3\tenex 
  \textfont\itfam=\nineit\def\it{\fam\itfam\nineit}%
  \textfont\slfam=\ninesl\def\sl{\fam\slfam\ninesl}%
  \textfont\ttfam=\ninett\def\tt{\fam\ttfam\ninett}%
  \textfont\gothfam=\ninegoth\scriptfont\gothfam=\sixgoth 
  \scriptscriptfont\gothfam=\fivegoth
  \def\goth{\fam\gothfam\tengoth}
  \textfont\bffam=\ninebf\scriptfont\bffam=\sixbf
  \scriptscriptfont\bffam=\fivebf
  \def\bf{\fam\bffam\ninebf}%
  \tt\ttglue=.5em plus.25em minus.15em
  \normalbaselineskip=11pt \setbox\strutbox\hbox{\vrule
  height8pt depth3pt width0pt}%
  \let\big=\ninebig\normalbaselines\rm}

\def\eightpoint{\def\rm{\fam0\eightrm}
  \textfont0=\eightrm \scriptfont0=\sixrm
  \scriptscriptfont0\fiverm
  \textfont1=\eighti \scriptfont1=\sixi
  \scriptscriptfont1\fivei 
  \textfont2=\eightsy \scriptfont2=\sixsy
  \scriptscriptfont2\fivesy 
  \textfont3=\tenex \scriptfont3=\tenex
  \scriptscriptfont3\tenex 
  \textfont\itfam=\eightit\def\it{\fam\itfam\eightit}%
  \textfont\slfam=\eightsl\def\sl{\fam\slfam\eightsl}%
  \textfont\ttfam=\eighttt\def\tt{\fam\ttfam\eighttt}%
  \textfont\gothfam=\eightgoth\scriptfont\gothfam=\sixgoth 
  \scriptscriptfont\gothfam=\fivegoth
  \def\goth{\fam\gothfam\tengoth}
  \textfont\bffam=\eightbf\scriptfont\bffam=\sixbf
  \scriptscriptfont\bffam=\fivebf
  \def\bf{\fam\bffam\eightbf}%
  \tt\ttglue=.5em plus.25em minus.15em
  \normalbaselineskip=9pt \setbox\strutbox\hbox{\vrule
  height7pt depth2pt width0pt}%
  \let\big=\eightbig\normalbaselines\rm}

\def\twelvepoint{\def\rm{\fam0\twelverm}
  \textfont0=\twelverm\scriptfont0=\tenrm
  \scriptscriptfont0\sevenrm
  \textfont1=\twelvei\scriptfont1=\teni
  \scriptscriptfont1\seveni 
  \textfont2=\twelvesy\scriptfont2=\tensy
  \scriptscriptfont2\sevensy 
   \textfont\itfam=\twelveit\def\it{\fam\itfam\twelveit}%
  \textfont\slfam=\twelvesl\def\sl{\fam\slfam\twelvesl}%
  \textfont\ttfam=\twelvett\def\tt{\fam\ttfam\twelvett}%
  \textfont\gothfam=\twelvegoth\scriptfont\gothfam=\ninegoth 
  \scriptscriptfont\gothfam=\sevengoth
  \def\goth{\fam\gothfam\twelvegoth}
  \textfont\bffam=\twelvebf\scriptfont\bffam=\ninebf
  \scriptscriptfont\bffam=\sevenbf
  \def\bf{\fam\bffam\twelvebf}%
  \tt\ttglue=.5em plus.25em minus.15em
  \normalbaselineskip=12pt \setbox\strutbox\hbox{\vrule
  height9pt depth4pt width0pt}%
  \let\big=\twelvebig\normalbaselines\rm}

\font \ninerm=cmr9
\font \medium=cmbx10 scaled \magstep1

\font \big=cmbx10 scaled \magstep2

\def\F{{\rm I\kern-.185em F}}
\def\R{{\rm I\kern-.185em R}}
\def\P{{\rm I\kern-.185em P}}
\def\N{{\rm I\kern-.185em N}}
\def\Z{{\mathchoice{\hbox{ Z\kern-.38em Z}}{\hbox{ Z\kern-.37em Z}}
{\hbox {\sevenrm Z\kern-.39em Z}}{\hbox{\sevenrm Z\kern-.39em Z}}}}
\def\Q{{\mathchoice{\hbox{\rm\kern.37em\vrule height1.4ex width.05em 
depth-.011em\kern-.37em Q}}{\hbox{\rm\kern.37em\vrule height1.4ex width.05em 
depth-.011em\kern-.37em Q}}{\hbox{\sevenrm\kern.37em\vrule height1.3ex 
width.05em depth-.02em\kern-.3em Q}}{\hbox{\sevenrm\kern.37em\vrule height1.3ex
width.05em depth-.02em\kern-.3em Q}}}}
\def\C{{\mathchoice{\hbox{\rm\kern.37em\vrule height1.4ex width.05em 
depth-.011em\kern-.37em C}}{\hbox{\rm\kern.37em\vrule height1.4ex width.05em 
depth-.011em\kern-.37em C}}{\hbox{\sevenrm\kern.37em\vrule height1.3ex 
width.05em depth-.02em\kern-.3em C}}{\hbox{\sevenrm\kern.37em\vrule height1.3ex
width.05em depth-.02em\kern-.3em C}}}}

\def\k{{\goth k}}
\def\b{{\goth b}}
\def\c{{\goth c}}
\def\h{{\goth h}}
\def\g{{\goth g}}
\def\t{{\goth t}}
\def\l{{\goth l}}\def\f{{\goth f}}
\def\m{{\goth m}}

\def\o{{\goth o}}
\def\p{{\goth p}}
\def\q{{\goth q}}

\def\s{{\goth s}}

\def\uu{{\goth u}}
\def\a{{\goth a}}
\def\X{{\bf X}}

\def\G{G^{\C}}

\def\xx{{\bf x}}\def\nn{{\bf n}}\def\mm{{\bf m}}
\def\yy{{\bf y}}

\def\RRe{{\rm Re}}
\def\RIm{{\rm Im}}

\def\pn{\par\noindent}
\def\sn{\smallskip\noindent}
\def\mn{\medskip\noindent}
\def\bn{\bigskip\noindent}
\def\pf{{\bf Proof. }}
\def\rem #1{}

\def\today{Roma, den\space\number\day\space\ifcase\month\or Januar 
\or Februar\or
M\"arz\or April\or Mai\or Juni\or Juli\or August\or
September\or Oktober\or November\or Dezember\fi,
\number\year}

\vglue 2cm

\pn
\centerline{{\bf \big{Complex extensions of semisimple symmetric spaces.}}}
\bn
\centerline{ \medium{ Laura Geatti}}
\centerline{ { Dipartimento di Matematica, Universit\`a di Roma 2 Tor Vergata, 00133 Roma }}
\centerline{ { e-mail: geatti@mat.uniroma2.it }}
\bn\bn\bn 
{\bf Abstract:} {\eightpoint Let $G/H$ be a pseudo-Riemannian semisimple symmetric space. The tangent bundle $T(G/H)$ contains a  
maximal $G$-invariant neighbourhood $\Omega$ of the zero section where the adapted-complex 
structure exists. Such $\Omega$ is endowed with a canonical $G$-invariant pseudo-K\"ahler metric of the same signature as the metric on 
$G/H$. 
We use  the polar map $\phi\colon \Omega \rightarrow G^\C/H^\C$  to define a $G$-invariant pseudo-K\"ahler metric
on  distinguished $G$-invariant domains in $G^\C/H^\C$ or on coverings of principal orbit strata in $G^\C/H^\C$. 
In the rank-one case, we show that the polar map is  globally injective and  the 
domain $\phi(\Omega)\subset G^\C/H^\C$ is an increasing union of $q$-complete domains.}

\beginsection Introduction.

Let $G/H$ be a non-compact semisimple symmetric space embedded in its complexification $G^\C/H^\C$. It is natural to ask whether there 
exists a $G$-invariant open set  
$$G/H~\subset ~D~\subset ~G^\C/H^\C,$$
whose complex analytic properties reflect the geometry and the harmonic analysis of $G/H$.
Since   $G^\C/H^\C$ contains $G$-invariant open sets with  very different complex analytic properties (cf.[Ge]), one should expect $D$ 
to be a proper subdomain of $G^\C/H^\C$.
The situation is best understood in the case of
 an irreducible   Riemannian symmetric space  $G/K$.  In this case   there is a distinguished  $G$-invariant subdomain  
$$G/K~ \subset ~D~\subset ~G^\C/K^\C,$$ which in many respects  may be considered the canonical complexification of $G/K$. 
The domain $D$, introduced by Akhiezer, Gindikin in [AG] and  intensively studied in recent years, has several remarkable properties: 
every eigenfunction of the algebra of $G$-invariant differential operators on $G/K$ admits a holomorphic extension to $D$ and every 
unitary
spherical representation of $G/K$ can be realized on a Hilbert space of holomorphic functions on $D$.  Moreover  $D$  is Stein,  
carries plenty of $G$-invariant plurisuharmonic functions, and is a maximal connected set where the $G$-action is proper (see 
[B],[BHH],[H],[KS1],[KS2]).  

 The Akhiezer--Gindikin domain $D$ can be described as follows. 
Let $\g$ be the Lie algebra of $G$ and let $\g=\k\oplus \p$ be a Cartan decomposition. The tangent bundle $T(G/K)$ of $G/K$ may be 
identified with the homogeneous vector bundle
$  G\times_K\p\rightarrow G/K$, where $K$ acts on $\p$ by the Adjoint representation. Then the  polar map 
$$\phi \colon G\times_K\p\longrightarrow G^\C/K^\C,\quad [g,X]\mapsto g\exp iX K^\C$$
may be viewed as a $G$-equivariant map from $T(G/K)$ with values in $G^\C/K^\C$.
The domain  $D$ is the image in $G^\C/K^\C$ of the largest connected neighbourhood $\Omega$ of $G/K$  in $T(G/K)$, where the  
differential of the  polar map $\phi$  has maximum rank. 
It turns out that the domain $\Omega$ also coincides with the maximal connected neighbourhood of 
$G/K$ in $T(G/K)$ where the so-called adapted complex structure  exists (see [LS],[Sz1] [GS1], [GS2]).  As a result,  $\Omega$ carries 
a canonical $G$-invariant 
K\"ahler structure  extending the Riemannian structure of $G/K$. Since the  polar map is globally injective on $\Omega$,  such a 
 $G$-invariant K\"ahler structure can be pushed-forward onto $D$. This makes the domain $D$ interesting also  from the  geometric point 
of view.

In this paper we consider  a semisimple pseudo-Riemannian symmetric space $G/H$, embedded in its complexification $\G/H^\C$, 
with the aim of determining how  the above facts generalize to this situation.  
If $G/H$ is a compactly causal symmetric space, then $G^\C/H^\C$ contains a maximal $G$-invariant Stein domain $D$ having $G/H$ in its Shilov boundary. The domain $D$, which may be considered as a one-sided complexification of $G/H$, has been investigated in [Ne1]. 
 
 However, in the general case the Stein manifold $\G/H^\C$ might contain no Stein $G$-invariant subdomains,
and that the $G$-action might fail  to be  proper on every invariant open subset of $G^\C/H^\C$. This happens for example  when $G/H$ is a real hyperboloid $SO(p,q)_0/SO(p-1,q)$, with $p,q>2$.  
Let $\g=\h\oplus \q$ be the decomposition of $\g$ induced by the symmetry of $G/H$ at the base point.  Identify the tangent bundle 
$T(G/H)$ with the homogeneous vector bundle $ G\times_H\q\rightarrow G/H$, where $H$ acts on $\q$ by the Adjoint representation. 
Consider then  the corresponding polar map
$$\phi\colon   G\times_H\q  \longrightarrow G^\C/H^\C,\quad [g,X]\mapsto g\exp iX H^\C.$$
Like in the Riemannian  case, the maximal connected set 
$G/H\subset \Omega$  where the polar map  has maximum rank  is a proper  
$G$-invariant subdomain in $T(G/H)$. The domain $\Omega$  also coincides with the maximal connected neighbourhood of $G/H$ in $T(G/H)$ 
where the adapted complex 
structure exists (cf. [Sz2], [HI]). As a result, $\Omega$ carries a canonical $G$-invariant pseudo-K\"ahler structure extending the 
pseudo-Riemannian structure 
of $G/H$.

If $G/H$  is a pseudo-Riemannian symmetric space of rank one, we show that the restriction of the polar map 
 $\phi\colon \Omega\longrightarrow G^\C/H^\C$  is globally injective.  
In particular, 
the  canonical $G$-invariant pseudo-K\"ahler structure of $\Omega$ can be pushed forward onto $D=\phi(\Omega)$. Moreover, if the metric 
on $G/H$ has signature $(p,q)$,  
the domain $D$ an increasing union 
of  $q$-complete domains (by our convention, $0$-complete is Stein).  In general,  $D$ is not a Stein domain.

In the higher rank case, the polar map $\phi\colon \Omega\longrightarrow G^\C/H^\C$ in generally not injective.  However,  $\phi$ is 
injective on close  orbits of maximal dimension in $\Omega$. Moreover, the restriction of $\phi$ to 
distinguished $G$-invariant subsets of $\Omega$ defines $G$-equivariant coverings of principal orbit strata in $D=\phi(\Omega)$. 

 As a result,  a canonical $G$-invariant pseudo-K\"ahler structure is defined on coverings of principal 
orbit strata in $D$ or on suitable neighbourhoods of closed orbits of maximal dimension in $D$. These results extend the ones 
obtained by Fels in the group case, using different methods (cf. [Fe]).

\sn
The paper is organized as follows. In section 1, we recall some definitions and  set up the notation.  
In section 2, we give several characterizations of the singular set of the differential of the polar map and we define the 
distinguished  
$G$-invariant neighbourhood $\Omega$ of $G/H$ in its tangent bundle $T(G/H)$. In section 3, we briefly recall the definition and the 
main properties of 
the adapted complex structure on $T(G/H)$.
In section 4, we prove two preliminary lemmas about the $G $-action on $\G/H^\C$, which may be of independent interest.  These lemmas are used to prove the main results in sections 
5 and 6.   In section 5, we deal with semisimple symmetric spaces of rank one.  We also work out in detail a family of examples. In section 6, we deal 
with symmetric spaces of  rank higher than one.

\beginsection 1. Preliminaries.

A semisimple symmetric space is a coset space $G/H$, where $G$ is a real
semisimple Lie group and $H\subset G$ is an open subgroup of the fixed point group of an involution $\tau\colon G\longrightarrow G$.

In what follows, we  consider {\it semisimple symmetric spaces  $G/H$ which admit a $G$-equivariant embedding into a simply connected 
complexification $G^\C/H^\C$}. They arise in the following way.
Start with a  simply connected complex semisimple Lie group $G^\C$ endowed with a Cartan involution $\Theta$, a conjugation $\sigma$ (different from  
$\Theta$) and a holomorphic involution $\tau$ satisfying the commutativity relations
$$\sigma\tau=\tau\sigma,\quad \Theta\sigma=\sigma\Theta,\quad \Theta\tau=\tau\Theta.\eqno(1.1)$$
Denote by $U=Fix(\Theta,\G)$ the corresponding compact real form, by $G=Fix(\sigma,\G)$ the corresponding non-compact real form  and by 
$H^\C=Fix(\tau,\G)$ the complex fixed point subgroup of $\tau$. By (1.1), the restriction of
$ \Theta$ to $G$ defines a Cartan involution $\theta$ of $G$ so that the maximal compact subgroup of $G$ is given by $K=G\cap U$. Similarly, 
the restriction of  $\tau$ to $G$ defines an involution  of $G$ commuting with $\theta$, whose fixed point subgroup is given by $H=G\cap 
H^\C$. 
In this way, the space $G/H$  admits an equivariant embedding in the complex symmetric space $G^\C/H^\C$  as the $G$-orbit of the base 
point $eH^\C$.

The product involution $\sigma^c:=\sigma\tau$ defines a conjugation of $G^\C$ with real form denoted by $G^c$; 
since $\sigma^c\tau=\tau\sigma^c$,   the restriction of
$\tau$ to $ G^c$ defines an involution of $G^c$, with fixed point subgroup $G^c\cap H^\C=H$. The restriction of
$\Theta$ to $G^c $ defines a Cartan involution  $\theta^c$ of $G^c$, commuting with $\tau$.
The  $G$-orbit and $G^c$-orbit  of the base point $eH^\C\in \G/H^\C$ define transversal totally real embeddings
$$G/H\hookrightarrow \G/H^\C \hookleftarrow G^c/H $$
of  so-called {\it c-dual symmetric spaces} [HO]. One has that
$$\dim_\R G/H=\dim_\R G^c/H=\dim_\C G^\C/H^\C.$$

Troughout the paper, the Lie algebra of a group is denoted by the corresponding gothic letter.
For example, $\g$ and $\g^\C$ denote the Lie algebras of $G$ and $\G,$ respectively.
An involution of a group  and the derived involution of its Lie algebra are denoted  by the same symbol. 
The commutativity relations (1.1) ensure  that the decompositions induced by $\Theta$, $\sigma$ and $\tau$ on $\g^\C$ and by their 
restrictions on $\g$ are all compatible with each other.
For example, if  $\g=\k\oplus \p$ is the Cartan
decomposition of $\g$ and $\g=\h\oplus \q$ is the decomposition of $\g$ induced by $\tau$, then 
both $\h$ and $\q$ are $\theta$-stable and $\g$ has a combined decomposition
$$\g=\h\oplus \q=\k\oplus \p=\h\cap \k~\oplus ~\h\cap \p~\oplus ~\q\cap \k~\oplus~ \q\cap \p.$$
In this setting, the $c$-dual symmetric spaces $G/H$ and $G^c/H$ are the analogues of  the non-compact Riemannian symmetric 
space $G/K$ and its compact dual $U/K$ in $\G/K^\C$.
The decompositions of $\g$ and $\g^c$ by $\tau$ are given  by
$$\g=\h\oplus\q,\quad\hbox{and} \quad \g^c=\h\oplus i\q,$$
respectively.
If $(\g=\h\oplus \q, \tau)$ is a symmetric algebra, a Cartan subspace of $\q$ is by definition a maximal abelian subspace $\c\subset 
\q$ 
consisting of semisimple elements.
The rank of a symmetric space $G/H$ is the dimension of an arbitrary Cartan subspace in $\q$.
\sn
In this paper, we do not deal with the group case, i.e. with symmetric spaces of the form
$G\times G/Diag(G)$.  Such spaces have already been investigated in [Br2],[Fe].

\beginsection 2. A distinguished $G$-invariant neighbourhood of $G/H$ in $T(G/H)$.

Let $G/H$ be a semisimple symmetric space. 
Identify the tangent bundle $T(G/H)$ with the homogeneous vector bundle 
$G\times_H\q$, defined as the quotient of  $ G\times \q $  by the $H$-action $h\cdot(g,X):=(gh^{-1}, Ad_hX)$. 
Identify $G/H$ with the zero section in $G\times_H\q$. In this way,   the {\it polar map}
$$ \phi\colon G\times_H \q \longrightarrow \G/H^\C,\quad [g,X]\mapsto g\exp iX  H^\C, \eqno(2.1)$$
 defines  a $G$-equivariant   map from the tangent bundle of $G/H$ with values in $ \G/H^\C$.
 \bn
The next proposition gives several characterizations of the set where the map $\phi$
has non-singular differential. By the $G$-equivariance of $\phi$ it is sufficient
to consider the differential $d\phi$ at the points $[e,X]$, with $X\in \q$.

\proclaim Proposition 2.1. Let $\phi$ be the map defined in (2.1). 
\item{(i)} Let  $X\in\q $. Then the differential $d\phi_{[e,X]}$ is non-singular if and only if
$$Ad_{\exp iX}\h\cap i\q=\{0\}.\eqno(2.2) $$
\item{(ii)} Let  $X\in\q $. Then $Ad_{\exp iX}\h\cap i\q=\{0\} $ if and only if $ ad_X\colon \g\rightarrow \g $ has no
real eigenvalue
  $$ \lambda\in \R,\quad \lambda\equiv\pi/2~{\rm mod}~\pi.$$
\item{(iii)} Let $X\in\q$ and let  $X=X_s+X_n $ be its Jordan decomposition, with $X_s$ semisimple,  $X_n$ nilpotent, $X_s$, 
$X_n\in\q$. 
Then 
 $Ad_{\exp iX}\h\cap i\q=\{0\}$ if and only if  $$ Ad_{\exp iX_s}\h\cap i\q=\{0\}.\eqno(2.3)$$ 
 If the semisimple element $X_s$ sits in a Cartan subspace $\c\subset \q$ and $\Delta_\c$ denotes the restricted root system of $\g^\C$ 
under $\c^\C$,  condition (2.3) is satisfied if and only if 
$$\alpha( X_s)\not\equiv\pi/2~{\rm mod}~\pi,\quad \hbox{for all}~ \alpha\in \Delta_\c.$$

 \mn\pf\pn
(i) The proof is similar to that of Prop.3 in [AG], where   the  compact real form $U$ of $G^\C$ is replaced 
by  the $c$-dual real form $G^c$. Let
$G^c/H$  the $c$-dual symmetric space of $G/H$. Observe that for $X\in \q,$ one has that
$$u=\exp iX\in G^c\quad{\rm and}\quad  u H\in G^c/H .$$
Consider the  diagram 
$$
\def\normalbaselines{\baselineskip20pt\lineskip3pt
\lineskiplimit3pt} 
\def\mapright#1{\smash{
\mathop{\longrightarrow}\limits^{#1}}}
\def\mapdown#1{\Big\downarrow
\rlap{$\vcenter{\hbox{$\scriptstyle#1$}}$}}
\matrix{G\times G^c &\mapright{\pi}&G\times_H (G^c/H)& \mapright{\phi}& G^\C/H^\C \cr 
\mapdown{1\times \rho_u^{-1}}&\quad&\quad& \nearrow\rlap{$\vcenter{\hbox{$\scriptstyle \tau_u$}}$}\cr
G\times G^c&\mapright{\psi_u}&G^\C/H^\C\cr}
$$
where the maps are defined as follows\pn
$\pi\colon G\times G^c\longrightarrow G\times_H G^c/H,\quad (g,v)\mapsto [g,vH]=[gh^{-1}, hvH],~~h\in H;$\pn
$\phi\colon G \times_H (G^c/H)\longrightarrow G^\C/H^\C,\quad [g,vH]\mapsto gvH^\C ;$\pn
$1\times \rho_{u^{-1}}\colon G\times G^c\longrightarrow G\times G^c,\quad (g,v)\mapsto (g,vu^{-1}),~~u\in G^c;$\pn
$\psi_u\colon G\times G^c\longrightarrow G^\C/H^\C,\quad (g,v)\mapsto   u^{-1} gvuH^\C,~~u\in G^c;$\pn
$\tau_u\colon G^\C/H^\C\longrightarrow G^\C/H^\C,\quad x\mapsto  ux H^\C,~~u\in G^c .$\pn
One easily checks  that the diagram is commutative:
$$\phi\pi(g,v)=\tau_u\psi_u 1\times \rho_{u^{-1}}(g,v),\qquad \hbox{for all}~ (g,v)\in G\times G^c,~\hbox{and}~ u\in G^c.$$
In order to determine the points $[e,X]\in\G\times_H\q$ where the differential $d\phi$ has maximum rank,
 we examine the rank of
$d\phi$ at the points $[e, uH],~u=\exp iX\in G^c.$ 
\pn
Since $1\times \rho_{u^{-1}}$ and $\tau_u$ are diffeomorphisms  and $d\pi$ is onto, such rank is maximum if and only if the  map 
$\psi_u$
 has differential of maximum rank at $(e,e)$. 
The differential $d\psi_{u,(e,e)}\colon \g\oplus \g^c\longrightarrow \q^\C$ is given by
$$d\psi_{u,(e,e)}(A,B)=Ad_{u^{-1}}(A+B)~~{\rm mod}~\h^\C,\quad (A,B)\in \g\oplus \g^c,\quad \g^c=\h\oplus  i\q,$$ and the kernel is 
given 
by
$$\ker d\psi_{u,(e,e)}=\{(C,-C)~|~C\in \h\}\oplus\{(0,Ad_uD)~|~D\in \h\}\oplus\{(E,0)~|~E\in  Ad_u i\h\cap\q \}.$$
The differential  $d\psi_{u,(e,e)}$
 has maximum rank if and only if $\dim (\ker d\psi_{u,(e,e)})=2\dim \h$. This happens precisely when $$Ad_u\h\cap i\q=\{0\}.$$

\sn
(ii) Let $X\in \q$. The operator $Ad_{\exp iX}\colon \g^c\longrightarrow \g^c$ can be written as
$$Ad_{\exp iX}=\exp ad_{iX}=\cos ad_X +i \sin ad_X.\eqno(2.4)$$
Since $\cos ad_X (\h)\subset \h$ and $i\sin ad_X (\h)\subset i\q$,
the condition $Ad_{\exp iX}\h\cap i\q=\{0\} $ is equivalent to the injectivity of the map $\cos ad_X |\h\colon\h\rightarrow\h$, namely
$$\cos ad_X (H)\not=0 , \qquad \hbox{for all}~ H\in \h , ~H\not=0.\eqno(2.5)$$
Indeed if $\cos ad_X  \colon\g\rightarrow\g$ is injective, condition (2.5) is clearly satisfied. Conversely, assume that 
$\cos ad_X  \colon\g\rightarrow\g$ is not injective.
Then  $\cos ad_X $ has an eigenvalue $\mu=0$ and
$ad_X$ has a {\it real} real  eigenvalue  $ \lambda\equiv \pi/2~\rm{mod~}\pi. $ 
In particular, there exists a $\lambda$-eigenvector  $Z\in \g$,   and  $\tau Z\not=\pm Z$. 
Since  $\ker \cos ad_X$   is $\tau$-stable, the vectors $H:=Z+\tau Z$ and $Q:=Z-\tau Z$ define non-zero elements in 
$\ker \cos ad_X \cap \h$ and $\ker \cos ad_X\cap \q$, respectively.
It follows that $\cos ad_X |\h\colon\h\rightarrow\h$ is not injective either. 
In conclusion, 
$\cos ad_X |\h$ is injective if and only if $\cos ad_X  \colon\g\rightarrow\g$ is injective, and this happens if and
only if $ad_X$ has no eigenvalues
$$\cases{  \lambda\in\R\cr  \lambda\equiv \pi/2~\rm{mod~}\pi.}$$
\sn
(iii) We already saw in (ii) that condition (2.2) is equivalent to the injectivity of the operator $\cos
ad_X \colon \g\longrightarrow\g$. So we need to prove that $\cos
ad_X $ is injective if and only if $\cos ad_{X_s} $ is injective. 
From the decomposition
 $X=X_s+X_n $ and the fact that $[X_s,X_n]=0,$ it follows that
$$\cos ad_{X}=\cos ad_{X_s} +(\cos ad_{X_s}(\cos ad_{X_n}-I)-\sin ad_{X_s}\sin ad_{X_n})\eqno(2.6).$$
Since $\cos ad_{X_s}$ is semisimple, $\cos ad_{X_s}(\cos ad_{X_n}-I)-\sin ad_{X_s}\sin ad_{X_n}$ 
is nilpotent and these operators commute,  equation (2.6) is the Jordan decomposition of $\cos ad_{X}$.
It follows that $\cos ad_X  $ is injective if and only if $\cos ad_{X_s}  $ is injective, as requested.

\bn 
{\bf Remark 2.2.}  By Proposition 2.1,  the regular set of $d\phi_{[e,X]}$ is a {\it proper $Ad_H$-invariant  subdomain}  of $\q$. A 
result by Halversheid (cf. [Ha], p.17), implies that the singular subset of   $d\phi_{[e,X]}$ in $\q$   disconnects $\q$.  
So the connected component of the regular set of $d\phi$  containing $G/H$   is a proper $G$-invariant subdomain 
 of $G\times_H\q$, namely
 $$\Omega=G\times_H \omega, \qquad \omega=\{ X\in \q~|~ |\lambda |< \pi/2,~ \hbox{for all}~ \lambda\in spec(ad_X)\cap \R\}.\eqno(2.7)$$
(Here $spec(L)$ denotes the spectrum of an operator $L$). Since $\omega$ is starlike, $\Omega$ is smoothly retractible to $G/H$.
\sn
By Proposition 2.1(ii), one has that $d\phi_{[e,0]}$ is non-singular
and  $d\phi_{[e,N]}$ is non-singular for every nilpotent element $N\in \q$.
In this framework, one can consider the  map $p\colon \q\longrightarrow \q ||Ad_H$, which associates to  $X\in\q$  the unique closed 
$Ad_H$-orbit in the closure of $Ad_H(X)$ (see [Br1]). Each fiber  of this map contains a unique closed orbit, which is also the unique 
orbit of minimum dimension. Recall that an $Ad_H$-orbit in $\q$ is closed if and only if $X$ is semisimple (cf. [vD]). By Proposition 
2.1(iii),  both $\omega$ and its boundary $\partial \omega$ are $Ad_H$-saturated sets, i.e. satisfy $p^{-1}p(\omega)=\omega$ (resp. 
$p^{-1}p(\partial\omega)=\partial\omega)$. 

 \proclaim Proposition 2.3.  Let $\q=\q\cap\k\oplus\q\cap\p$ be the Cartan decomposition of $\q$, let $\a\subset \q\cap\p$ be a maximal 
abelian subspace (not necessarily maximal abelian in $\q$) and let $\Delta_\a=\Delta_\a(\g,\a)$ be the corresponding restricted root 
system. Define $ \omega_0=\{A\in\a~|~|\alpha(A)|<\pi/2,~\forall \alpha\in\Delta_\a\}.$
  Then the set $\omega$ defined in (2.7) is given by
$$\omega=p^{-1}(p( Ad_H(\q\cap\k\oplus \omega_0)) ).$$

\mn\pf Both sets $Ad_H(\q\cap\k\oplus \omega_0)$ and $\omega$ are $Ad_H$-stable.
So we need to show that they have the same semisimple elements.
Let $X\in Ad_H(\q\cap\k\oplus \omega_0)$. If $\lambda\in spec(ad_X)\cap \R$, then $| \lambda |<\pi/2$. This shows that
$Ad_H(\q\cap\k\oplus \omega_0) \subset \omega$.
\pn
To prove the converse statement,  let $X\in\omega$ be a  semisimple element.  Then $X$ is $H$-conjugate to an element  
$S=S_\k+S_\p=Ad_hX$ 
in a standard Cartan subspace
$\c=\c_k\oplus\c_\p$, with $c_\k\subset \q\cap\k$ and $\c_\p\subset \a$ (see [Ma], p.79). In particular $S_\p\in \omega_0\subset \a$ 
and 
$X\in Ad_H(\q\cap\k\oplus \omega_0)$.

\beginsection 3. The adapted complex structure  on $\Omega$.

The notion of an adapted complex structure on the tangent bundle of a Riemannian manifold was introduced and developed
in [LS], [Sz1], [GS1], [GS2]. Its definition and many of its features have a straightforward generalization to the pseudo-Riemannian 
case 
[Sz2].
\bn
{\bf Definition 3.1}. Let $(M,g)$ be a pseudo-Riemannian manifold and let $U$ be an
open neighbourhood of $M$ in its tangent bundle $TM.$ A complex
structure $J$ on $U$ is called {\it adapted} if 
for every geodesic $\gamma\colon \R\longrightarrow M,$ the differential
$$d\gamma \colon T\R\cong \C \longrightarrow TM, 
\quad (x,y)\mapsto (\gamma(x), y\gamma'(x))$$ is holomorphic on 
$d\gamma^{-1}(U).$

\bn 
A semisimple symmetric space $G/H$ is in a natural way a pseudo-Riemannian manifold.
The tangent space $T(G/H)_{eH}$  to $G/H$ at the base point $eH$, can be identified with 
$\q=\q\cap \k \oplus \q\cap \p$;  the restriction of the Killing form of $\g$ to $\q\times \q$
$$B_\g|\q(X,Y)=Tr (ad_X\circ ad_Y),\quad X,Y\in\q$$ induces on $G/H$ a
$G$-invariant metric $g$ of signature $(\sigma^+,\sigma^-)$, where 
$\sigma^+=\dim \q\cap \p $ and $\sigma^-=\dim \q\cap \k.$  \sn

\bn
Let $\Omega$ be the domain in $  G\times_H\q\cong T(G/H)$ defined in (2.7). In analogy with the Riemannian case, one has that (cf. 
[HI], 
[Sz2])
\bn
{\sl \item{(a)} The complex structure on $\Omega$ given by the pull-back of the complex structure of $G^\C/H^\C$ by the polar map (2.1) 
is 
adapted. 
\item{(b)} $\Omega$ is the ``largest"  $G$-invariant connected subset of $T(G/H)$ containing $G/H$ and carrying an adapted complex 
structure.   
\item{(c)}  The energy function
$$E\colon T(G/H) \longrightarrow \R ,\quad  E(x,v):={1\over 2}g_x(v,v)^2,\quad x\in G/H,~v\in T(G/H)_x \eqno(3.1)$$ 
 is a smooth $G$-invariant function on  $\Omega$.  Its complex Hessian of $E$ has $\sigma^+$ positive and $\sigma^-$ negative 
eigenvalues. 
\item{(d)} The formula 
$$h(Z,W):=-{i\over 2}\partial\overline \partial E( Z,\overline W),\quad Z,W\in T^\C\Omega. \eqno(3.2)$$ 
defines a $G$-invariant pseudo-K\"ahler metric on $\Omega$ 
with the same signature as $g$.  
\item{(e)} The function $\sqrt{|E|}$ satisfies the homogeneous complex Monge-Ampere equation $$(\partial\bar 
\partial\sqrt{|E|})^n\equiv 
0,$$  outside the null set in $\Omega$.
\item{}}

\beginsection 4. The $G $-action on $\G/H^\C$.

In this section we prove two preliminary lemmas regarding the $G $-action on $\G/H^\C$. These are  used in the next sections to prove 
the 
main results.  

 \bn
Resume the notation introduced in section 1.  Denote by $Aut_\R(\g^\C)$ the group of the real automorphisms of the complex Lie algebra 
$\g^\C$. 
Define a  map $\eta\colon \G\rightarrow Aut_\R(\g^\C)$  by $\eta(x)=\sigma Ad_x\tau Ad_{x^{-1}}$  (see [Ma], p.51).  An element  $x\in 
\G$ 
is called   {\it regular semisimple with respect to $\sigma, \tau$} if its image 
 $\eta(x)$ is a regular semisimple element in the real algebraic group $Aut_\R(\g^\C)$ or equivalently when it  sits on a closed 
$G\times 
H^\C$-orbit of maximal dimension in $\G$.  Let  $\bar x$ denote  the image of  $x $ under the canonical projection
$\pi\colon G^\C\rightarrow G^\C/K^\C$.  
Then $x$ is regular  semisimple  with respect to $\sigma, \tau$  if and only if  its  image $\bar x\in \G/H^\C$ sits on a closed 
$G$-orbit 
of maximal dimension in $\G/H^\C$ (see [Ma]). 
Let  $X\in\q$, such that $x=\exp iX\in\G$ is a regular semisimple with respect to $\sigma, \tau$.
Then $X$ sits in some Cartan subspace $\c$ in $\q$ and can be characterized in terms of the restricted roots 
$\Delta_\c=\Delta_\c(\g^\C,\c^\C)$ as follows (see [Ge], Prop. 3.14): $X$ sits in the complement in $\c$ of the  set
$$ \bigcup_{\alpha\in \Delta_\c^r}\{ \alpha(X) \equiv 0~{\rm mod}~\pi/2\} 
\bigcup_{\alpha\in \Delta_\c^i}\{\alpha(X) =0 \}\bigcup_{\alpha\in \Delta_\c^c}
\Bigl\{\cases{\RRe\alpha(X) \equiv 0~{\rm mod}~\pi/2 \cr \RIm\alpha(X)=0 } \Bigr \}\eqno(4.1)$$
(here $\Delta_\c^r$, $\Delta_\c^i$, $\Delta_\c^c$ denote the sets of roots which restricted to $\c$ take real, imaginary or complex 
values, respectively). 
Observe that one such $X$ is in particular regular semisimple in $\q$, which by definition means that
the centralizer $Z_\q(X) $ is abelian and  equal to $\c$. 
This last condition is characterized by
$\alpha(X)\not= 0$, for all $\alpha\in \Delta_\c$.
\bn
{\bf Remark 4.1.}  
Let $X=X_s+X_n$  be the Jordan decomposition of $X$  in $\q$. Let $x=\exp iX= x_sx_n$  be the corresponding Jordan decomposition in 
$\G$, 
with $x_s=\exp iX_s$ and $x_n=\exp i X_n$.  
  One has that 
 $$\eta(x)=Ad_{x_n^{-2}}\sigma\tau Ad_{x_s^{-2}}=Ad_{x_n^{-2}}\eta(x_s).$$
 Set $u=Ad_{x_n^{-2}}$ and $s=\sigma\tau Ad_{x_s^{-2}}=Ad_{x_s^{-2}}\sigma\tau.$ 
It is easy to check that
$u$ is unipotent, $s$ is semisimple, and $su=us$. So  $\eta(x)=us$ is the Jordan decomposition of $\eta(x)$ in $Aut_\R(\g^\C)$. 
Moreover, 
one has that $$\sigma (iX_n)=Ad_{x_s}\tau Ad_{x_s^{-1}}(iX_n)=- iX_n.$$ 
Then by [Ma], Prop.2 (ii), p.66, the decomposition  $x=x_sx_n$  also coincides with the lifting to $\G$ of the Jordan decomposition  of 
$\eta(x)$   in $Aut_\R(\g^\C)$.

\proclaim Lemma 4.2. Let $X \in \q$.   Let $x=\exp iX$ be the corresponding element  in $\G$ and $\bar x $ its image  in $\G/H^\C$. 
Then
\item{(i)} The isotropy subgroup of $\bar x$ in $G$ is given by $G_{\bar x}=\{g\in G~|~ gx^2=x^2\tau(g)\} $;
\item{(ii)} $G_{\bar x}\subset  Z_G(x^4)$;
\item{(iii)} $Z_G(x^4)=Z_G(x^2)$ implies $G_{\bar x}=Z_H(x^2)  $;
\item{(iv)} Assume that $[e,X]\in G\times_H\q$ lies in the regular set of $d\phi$.  Then 
  $G_{\bar x}=Z_H(x^2) $;
\item{(v)}  Let  $\c$ be  a  Cartan subspace in $\q$ and let $X\in\c$. Assume that $x=\exp iX$ is a regular semisimple
element with respect to
$\sigma,\tau$. Then 
 $G_{\bar x}=Z_H(x )= Z_H(X)=Z_H(\c).$ 

\mn
\pf\pn
(i) By definition, $g\in  G_{\bar x}$ if there exists $h_c\in H^\C$ such that 
$gx=x h_c$.
Write $h_c=x^{-1}g x$. Since $h_c\in H^\C$, one has that $\tau(h_c)=h_c$. This is equivalent to 
$$  x\tau(g) x^{-1}=x^{-1}gx \quad \hbox{and} \quad x^2\tau(g)=gx^2.$$
Conversely, assume that $g\in  \{g\in G~|~ gx^2=x^2\tau(g)\}$. Then 
$$gx=gx^2 \cdot x^{-1}=x^2\tau(g)x^{-1}=x\cdot x\tau(g) x^{-1}.$$
Define $h_c:=x \tau(g)x^{-1}$. One has that 
$$\tau(h_c)=\tau(x \tau(g)x^{-1})=x^{-1}gx=x^{-1}gx^2x^{-1}=x^{-1}x^2\tau(g)x^{-1}=x\tau(g)x^{-1}=h_c.$$
So  $h_c\in H^\C$ and $g\in  G_{\bar x}$, as desired.
\pn
(ii) Let $g\in G_{\bar x}$. Then by (i) one has  that $gx^2=x^2\tau(g)$ and $g=x^2\tau(g)x^{-2}$. Since $g\in G$,  one has that 
$\sigma(g)=g$, which is equivalent to
$$  \sigma(x^2)\tau(g)\sigma(x^{-2})=  x^2 \tau(g) x^{-2} 
\quad \Leftrightarrow \quad x^4\tau(g)=\tau(g)x^4\quad \Leftrightarrow \quad x^4g=gx^4.$$
In other words,  $g\in G_{\bar x}$ implies $g\in Z_G(x^4)$.
\pn
(iii) Assume that $Z_G(x^2)=Z_G(x^4)$. By (ii) one has that $G_{\bar x}\subset Z_G(x^2)$.
This together with (i) implies that $gx^2=x^2\tau(g)=x^2g$. Then $g\in H$ and $G_{\bar x}=Z_H(x^2)$.
 \pn
 (iv) By (iii) it is sufficient to show that  $Z_{\G}(x^4)= Z_{\G}(x^2)$, and actually 
 that $Z_{\G}(x^4)\subset  Z_{\G}(x^2)$, the opposite inclusion being obvious. Assume first that $X$ is a {\it semisimple} element in 
some 
Cartan subspace $\c\subset\q$. 
Since $\G$ is simply connected, the centralizers  $Z_{\G}(x^2)$ and $Z_{\G}(x^4)$  are connected (cf. [Hu]) 
and hence determined by their Lie algebras.  
 Denote by $\Delta_\c=\Delta(\g^\C,\c^\C)$ the   restricted root  system of $\g^\C$ with respect to $\c^\C$. Let 
 $\g^\C=Z_{\g^\C}(\c^\C)\oplus \bigoplus_{\alpha\in\Delta_\c}\g^\alpha $ be the corresponding  root decomposition.  The  Lie algebra of 
$Z_{\G}(x^4)$ is given by 
$$Z_{\g^\C}(x^4)=\{Z\in\g_\C~|~ Ad_{x^4}Z=Z\}=Z_{\g^\C}(\c) \oplus \bigoplus_{\alpha(X)=0}\g^\alpha\oplus 
\bigoplus_{\alpha(X)\not=0\atop
\alpha(4X)\equiv 0~{\rm
mod}~2\pi}\g^\alpha.$$
 By Prop. 2.1(iii), the element $[e,X]$ lies in the regular set of $d\phi$ if and only if  $\alpha(X)\not=(2k+1)\pi/2,~k\in\Z$, for all 
$\alpha\in\Delta_\c$. As a consequence,
$\beta(4X)\equiv 0~{\rm mod}~2\pi$, for some root $\beta$, if and only if  
$$\beta(X)=m\pi/2,\quad\hbox{for some}~  m\in 2\Z,~m\not=0.$$
It follows that $\beta(2X)\equiv 0~{\rm mod}~2\pi$, which means that $Z_{\g^\C}(x^4) \subset
Z_{\g^\C}(x^2)$ and $Z_{\G}(x^4)\subset  Z_{\G}(x^2)$, as requested.

\mn Assume now that $X$ is {\it non-semisimple}. Let $X=X_s+X_n$ be its Jordan decomposition in $\q$, with $X_s,X_n\in\q$ and 
$[X_s,X_n]=0$.
Write  $x=x_sx_n=\exp iX_s\exp iX_n $. By Remark 4.1 and [Ma], Prop.2, p.66,  the equation
$$ gx_s x_n=x_s x_n h_c,\quad g\in G,~h_c\in H^\C$$
is equivalent to the system
$$\cases{ g\exp iX_s  =\exp iX_s   h_c\cr Ad_gX_n=X_n .}   $$
In particular,  $g\in G_{\bar x_s}\cap Z_G(X_n)$.
Since $[e,X_s]$ lies in the regular set of $d\phi$ (by Proposition 2.1(iii)),  one has that $g\in H$. It follows that
$$G_{\bar x}=G_{\bar x_s}\cap Z_G(X_n)= Z_H(x_s^2)\cap Z_G(X_n)\subset H,\quad\hbox{and}\quad  G_{\bar x}=Z_H(x^2),$$
as requested.
\pn
(v) If $x=\exp iX$ is a regular semisimple element with respect to $\sigma,\tau$, then $X $ is a semisimple element  in
some Cartan subspace
$\c\subset \q$ and $[e,X]$ lies in the regular set of
$d\phi$.   The same argument used in (iv) shows that 
$$  Z_{\g^\C}(x^4)=Z_{\g^\C}(\c^\C)  =Z_{\g^\C}(x )=Z_{\g^\C}(\c) =Z_{\g^\C}(X).$$
As a consequence, the $\sigma$-stable connected groups  
$Z_{\G}(x )$, $Z_{\G}(x^4)$, $Z_{\G}(\c)$ and  $Z_{\G}(X)$  all coincide (for the connectedness of such groups, see [Hu][St]). In 
particular (see also [Ge], Sect.3.3),  one has that 
$$G_{\bar x}=H\cap Z_{\G}(x )=Z_H(x )=H\cap Z_{\G}(X )= Z_H(X)=Z_H(\c).\eqno(4.2)$$

 \bn
\proclaim Lemma 4.3. Let $X,Y\in \q$. Let $x=\exp iX$, $y=\exp iY$  be the corresponding elements in $\G$ and $\bar x$, $\bar y$ their 
images in $\G/H^\C$.
One has that:
\item{(i)} $\bar x$, $\bar y$ sit on the same $G$-orbit if and only if   $gx^2=y^2\tau(g)$, for some $g\in G$;
\item{(ii)} $\bar x$, $\bar y$ sit on the same $H$-orbit in $\G/H^\C$ if in addition  $gx^2=y^2g$ holds; 
\item{(iii)} if $\bar x$, $\bar y$ sit on the same $G$-orbit,  then $y^4= gx^4g^{-1}$, for some $g\in G$.
\pn
Let $X,Y$ be semisimple elements in the same Cartan subspace $\c\subset \q$. Assume that  $x $, $y $ are regular semisimple elements 
with 
respect to $\sigma,\tau$ and that  $\bar x $, $\bar y $ sit on the same $G$-orbit, i.e. (4.3) holds.  Then:
\item{(iv)} $g\tau(g)^{-1}\in Z_H(\c)$ and $\tau(g)=zg$, for some $z\in Z_H(\c)$, with $z^2=1$;
\item{( v)} $Ad_gX\in\c$ and $g\in N_G(\c)$;  
\item{( vi)} $gxg^{-1}=yq$, with $g\in N_G(\c)$ and $q^4=1$.

\mn\pf\pn
(i) By definition, $\bar x$, $\bar y$ sit on the same $G$-orbit  in  $\G/H^\C$ if
$$gx=yh_c, \quad \hbox{ for some}~g\in G,~h_c\in H^\C.\eqno(4.3)$$
Write $h_c=y^{-1}gx $. Then $\tau(h_c)=h_c$ implies $gx^2=y^2\tau(g)$. Conversely, assume that $gx^2=y^2\tau(g)$, for some $g\in G$. 
Write 
$gx=gx^2  x^{-1}=y^2\tau(g)x^{-1}=y y\tau(g)x^{-1}$ and set $h_c=y\tau(g)x^{-1}$. One can check that $\tau(h_c)=h_c$  and $h_c\in 
H^\C$. 
This shows that (4.3) is satisfied and $\bar x$, $\bar y$ sit on the same $G$-orbit.
\pn
(ii)
If $g\in H$, then $\tau(g)=g$ and  one has $gx^2=y^2g$ by statement (i). Conversely, assume that (4.3) holds and that moreover 
$gx^2=y^2g$, for $g\in G$. Then by (i)  $\tau(g)=g$ and $g\in H$.
\pn
(iii) By (i), our assumption is equivalent to  $gx^2=y^2\tau(g)$, for some $g\in G$. Write  $g=y^2\tau(g)x^{-2} $. Then $\sigma(g)=g$ 
implies $\tau(g)x^4=y^4\tau(g)$. By applying the involution $\sigma\tau$
to both  terms of the equality, we get $y^4=Ad_gx^4$, as desired.
\pn
(iv) By applying the map $\eta_{\sigma\tau}(x):=x\sigma\tau(x)^{-1}$ to both terms of equation (4.3), we get
$$g\tau(g)^{-1}y=yh_c\sigma(h_c)^{-1}.$$
This means that $g\tau(g)^{-1}$ is an element in $G_{\bar y}$, the isotropy subgroup in $G$ of $\bar y$.  
Since $y $ is regular semisimple with respect to $\sigma,\tau$, by Lemma 4.2(v), one has that $g\tau(g)^{-1}\in Z_H(\c)$.  
Equivalently,  
$\tau(g)=zg$, for some $z\in
Z_H(\c)$. Since $\tau(z)=z=z^{-1}$, it follows that $z^2=1$.

\pn
(v) We first prove that $Ad_gX\in\q$. We need to show that $\tau(Ad_gX)=-Ad_gX$, which is equivalent to
 $g^{-1}\tau(g)\in Z_G(X)$.
From $\sigma\tau(y)=y$, $\sigma\tau(y^4)=y^4$ and  $y^4=Ad_gx^4$ (by (iii)), we obtain   
$$ y^4=Ad_gx^4 = Ad_{\tau(g)}x^4 \quad \hbox{and} \quad g^{-1}\tau(g)\in Z_G(x^4).$$
Since $x $ is regular semisimple with respect to $\sigma,\tau$, by Lemma 4.2(v), one has that   
 $Z_G(x^4)=Z_G(X)$.  Then $g^{-1}\tau(g)\in Z_G(X)$ and  $Ad_gX\in \q$, as requested.
 \pn
 By (4.1), the operator $ad_{4Y}$ has no real eigenvalue $\lambda=k2\pi,~k\in\Z \setminus\{0\}$. 
 As a consequence,  
 the element $i4Y$ lies in the regular set of the differential of the exponential map $\exp\colon \g^\C\rightarrow \G $ 
(cf.[Va]). From
 $\exp i4 Y=Ad_g \exp i4 X=\exp i4 Ad_gX $, by  [Ne2],  one obtains   
$$[Y,Ad_gX]=0, \quad \hbox{or equivalently}\quad Ad_gX\in Z_\q(Y)=\c .$$
Since $X$, $Ad_gX$ are   regular semisimple in $\q$, from
$$\c=Z_\q(X)= Z_\q(Ad_gX)=Ad_gZ_\q(X)=Ad_g\c,$$
it follows that $g\in N_G(\c)$, as requested.
\pn
(vi) By (i) and (iv) we can write $gx^2g^{-1}=y^2\tau(g)g^{-1}=y^2z$, for some $z\in Z_H(\c)$, with $z^2=1$. By (v), we can also write 
$\exp i2(Ad_gX-Y)=\exp \gamma$, for some $\gamma$ in the square  lattice ${1\over 2}\Gamma$ (here $\Gamma$ denotes the unit lattice in 
$i\c\cap \uu\subset \g^\C$). 
It follows that $Ad_gX=Y+\eta$, for some $\eta \in {1\over 4}\Gamma$, and 
$gxg^{-1}=yq$, for some $q\in \exp i\c$, with $q^4=1$.

\beginsection 5.1. The rank-1 case.  

Let $G/H$ be a semisimple symmetric space of rank one. Let $\Omega\subset G\times_H\q$ be the domain  defined in (2.7). The main goal 
of 
this section is to prove that the polar map $\phi\colon \Omega\rightarrow \G/H^\C$ is globally injective (cf. Proposition 5.4). This 
shows 
that the image domain $D=\phi(\Omega)$ in $G^\C/H^\C$ is the direct generalizazion of the Akhiezer-Gindikin domain. However, in the 
pseudo-Riemannian case the domain $D$ is generally not Stein. 
\bn
Recall that in tthe rank-one case every Cartan subspace 
in $ \q$ is one-dimensional.  Up to $Ad_H$-conjugacy, there are precisely two $\theta$-stable Cartan subspaces
in $\q$: a compact one $\t\subset \q\cap\k$ and a non-compact one $\a\subset \q\cap\p$.
As a consequence, given  a semisimple element in $X\in \q$, the eigenvalues of
$ad_{ X}$  are either all real or all imaginary. 
Let $\a  $ be the non-compact Cartan subspace in $\q$  and let $\Delta_\a=\Delta_\a(\g,\a)$ denote the corresponding restricted root 
system.
Then either  $\Delta_\a=\{\pm \alpha\}$ is  of type $A_1$ (reduced case) or  $\Delta_\a=\{\pm \alpha,\pm 2\alpha\}$ is of typre $BC_1$ 
(non-reduced case). A list of all rank-1 semisimple symmetric algebras and their restricted root systems can be found in [OS].
Here is the list:
\bn
$\bullet$ $\s\o(p+1 ,q+1 )_0/\s\o(p+1 ,q ) $, for $p\ge 0,~q\ge 0$
\pn
$q=0$ Riemannian
\pn
$\Delta_\a=A_1$.

\bn
$\bullet$ $\s\uu(p+1,q+1)/\s(\uu(p+1,q )\times \uu(1)) $, for $p\ge 0,~q\ge 0$
\pn
$q=0$ Riemannian
\pn
$\Delta_\a=BC_1$.

\bn
$\bullet$  $\s\p(p+1,q+1)/\s\p(p+1,q )\times \s\p(1)) $, for $p\ge 0,~q\ge 0$
\pn
$q=0$ Riemannian
\pn
$\Delta_\a=BC_1$.

\bn
$\bullet$ $\s\l(n+2,\R)/\g\l(n+1,\R)$, for $n\ge 0$
\pn
$\Delta_\a=BC_1$.

\bn
$\bullet$ $\s\p(n+2,\R)/\s\p(n+1,\R)\times \s\p(1,\R)$, for $n\ge 0$
\pn
$\Delta_\a=BC_1$.

\bn
$\bullet$ $\f_{4(-20)}/\s\o(8,1)$
\pn
$\Delta_\a=BC_1$.

\bn
$\bullet$ $\f_{4(4)}/\s\o(5,4)$
\pn
$\Delta_\a=BC_1$.

\bn 
Let $\l =\b\oplus\a=\b_\k\oplus\b_\p\oplus \a $ be a  $\theta$-stable Cartan subalgebra of $\g$ extending $\a$ and let
$\Delta=\Delta(\g^\C,\l^\C)$ denote the corresponding root system. By [OS], Sect. 3.2,  there exists a positive system $\Delta^+$ in 
$\Delta$ which is stable both  under $-\tau$ and $-\theta$.  Let $\Pi=\{\lambda_i\}_{i=1,\ldots,n}$ be the corresponding set of simple 
roots. 
Let $\Gamma$ be the inverse  lattice  in $i\l_\R$
$$\Gamma=\bigoplus_{\lambda_i\in \Pi}\Z2 \pi ih_{\lambda_i}\subset i\l_\R, \quad \lambda_i(h_{\lambda_i})=2.$$
Recall that in a simply connected compact Lie group $U$, the lattice $\Gamma$ coincides with the unit lattice $\exp_U^{-1}(e)$. Denote 
by 
$\langle \cdot,\cdot\rangle$ the euclidean inner product induced on $\l_\R$ by the restriction of the Killing form of $\g$. Then  the 
following relations hold
$$\Vert h_{\lambda_i}\Vert={4\over {\Vert\lambda_i \Vert}},\quad 
\langle h_{\lambda_i},h_{\lambda_j}\rangle ={2\over {\Vert\lambda_i \Vert}}c_{ij}={2\over {\Vert\lambda_j \Vert}}c_{ji},\qquad 
c_{ij}\in\{0,-1,-2,-3\}  .$$

\bn
\proclaim Lemma 5.1. Let $\Delta$ be a root system, all of whose roots have the same lenght (namely $A_n,~D_n,~E_6~,E_7,~E_8$).
\item{(i)} The  generators $\{i \pi h_{\lambda_i}\}$  of ${1\over 2}\Gamma$ are vectors of shortest length in ${1\over 2}\Gamma$.
\item{(ii)} Let $\tau$ be an involutive automorphism of ${1\over 2}\Gamma$, which is  an isometry  and satisfies the conditions
$$\langle h_{\lambda_i},\tau h_{\lambda_i}\rangle\ge 0,\quad i=1,\ldots,n.$$ Then 
$X:=i\pi (h_{\lambda_i}-\tau h_{\lambda_i})$ is the shortest vector in $ \R X\cap {1\over 2}\Gamma$.

\mn\pf  
Without loss of generality we can prove the corresponding statements for  a lattice $L$
with basis vectors $\{v_1,\ldots, v_n\}$, satisfying
$$\Vert v_i\Vert^2=2,\qquad \langle v_i,v_j\rangle=0,-1.$$
\item{(i)}  Let $v=\sum_im_iv_i,~m_i\in\Z$. The  formula 
$$  \Vert v\Vert^2=2(\sum_im_i^2-\sum_{i<j}\epsilon(i,j) m_im_j)\ge 2 ,\quad \epsilon(i,j)=0,1 $$
show that the square of the norm of every vector in the lattice satisfies is a positive even number.  Hnece the basis vectors are 
vectors 
of shortest length as requested.
\item{(ii)} Write  $v:= v_i-\tau v_i$. The formula
$$\Vert v\Vert^2= \Vert v_i\Vert^2+ \Vert \tau v_i\Vert^2-2\langle v_i,\tau v_i\rangle=2( \Vert  v_i\Vert^2-\langle v_i,\tau 
v_i\rangle)  
0 $$
implies that the inner product
 $\langle v_i,\tau v_i\rangle $ can only take the values 1 or 0.  In the first case, $\Vert v\Vert^2=2 $ and $v$ is a vector of 
shortest 
length. In the second case, $\Vert v\Vert^2=4$. Assume that there is a  vector $v'\in \R v\cap L$ which is shorter than $v$. By (i), 
$\Vert v'\Vert^2 =2$.   Since $v'=q v$, for some rational number $q\in \Q$, the equation
 $$\Vert v'\Vert^2 =q^24=2$$  implies that $ q={1\over \sqrt 2} $. This is absurd.

\bn
\proclaim Lemma 5.2. 
Let  $\Delta$ be the root system $B_n$. A  generator  $\{i\pi h_{\lambda_{i_0}}\}$  of $ {1\over 2}\Gamma$ is a vector  of shortest 
length 
if and only if $\lambda_{i_0}$ is a long root in $\Delta$.

\mn\pf Without loss of generality we can prove the corresponding statement  for  a lattice $L$ with basis
$\{v_1,\ldots,v_{n-1},v_n\}$ satisfying
$$ \Vert v_1\Vert^2=\ldots=\Vert v_{n-1}\Vert^2=2, ~\Vert v_n\Vert^2=4,$$
$$ \langle v_i,v_j\rangle=-1,~1\le i<j\le n-1,~~\langle v_{n-1},v_n\rangle=-2.$$
The first $n-1$ vectors correspond to the long roots in $\Delta$, while $v_n$ corresponds to the short root. Observe that
if $v=\sum_im_iv_i$ is an element in $ L$, then
$\Vert v\Vert^2$ is a positive  even integer. Consider then  $v=v_i-\tau v_i$, for some $i=1,\ldots,n-1$. The inequality
$$ \Vert v\Vert^2=2(\Vert v_i\Vert^2-\langle v_i,\tau v_i\rangle)>0,$$ 
implies that $\langle v_i,\tau v_i\rangle $ can only take the values 1 or 0. In the first case $\Vert v\Vert^2=2$ and $v$ is a vector 
of 
shortest length. In the second case, 
$\Vert v\Vert^2=4$. Since  no  rational number satisfies $q^2={1\over 2}$,  one has that $v=v_i-\tau v_i$ is the shortest vector in $  
\R 
v\cap L$.
On the other hand,  one easily checks that if  $v=v_n-\tau v_n$,   then  $ {1\over 2}v$ is the shortest vector in $ \R v\cap L$.

\bn\bn
\proclaim Lemma 5.3. Let $\gamma$ be a non-zero vector in $ \a\cap i{1\over 2}\Gamma$. Then for every root $\alpha$ in the restricted 
root 
system $\Delta_\a$, one has that
\item{(i)} $\alpha(\gamma)\in\Z \pi$;
\item{(ii)} $|\alpha(\gamma)|\ge \pi$;
\item{(iii)} $|\alpha(\gamma)|\ge 2 \pi$, if  $\Delta_\a$ is reduced. 

\mn\pf
It is sufficient to prove the lemma for the simple root in $ \alpha\in  \Delta_\a$. Recall that  $\alpha$ can be written as 
$\alpha={1\over 2}(\lambda_{i_0}-\tau\lambda_{i_0})$, for some simple root $\lambda_{i_0}\in \Pi$. Let $\gamma= \pi \sum m_i 
h_{\lambda_i}\in i{1\over 2}\Gamma$, with $m_i\in\Z$. Then
$$\alpha(\gamma)={1\over 2}(\lambda_{i_0}-\tau\lambda_{i_0})(\gamma)={1\over 2}(\lambda_{i_0}(\gamma)-\tau\lambda_{i_0}(\gamma))=
{1\over 2}(\lambda_{i_0}(\gamma)-\lambda_{i_0}(\tau\gamma))=$$
$$=\lambda_{i_0}(\gamma)=\pi \sum m_i  \lambda_{i_0}( h_{\lambda_i})\in  \pi\Z,$$
proving (i).
To prove (ii), we need to show that 
$$\alpha(\gamma)\not=0, \quad \hbox{for all}~ \gamma\in i{1\over 2}\Gamma\cap \a,~\gamma\not=0.$$
Observe that the vector  $\gamma_0= \pi(h_{\lambda_{i_0}}-\tau h_{\lambda_{i_0}} )$ is a non-zero vector   $ \a\cap i{1\over 2}\Gamma,$ 
and that  every $\gamma\in  \a\cap i{1\over 2}\Gamma$ is of the form $\gamma=q\gamma_0$, for some $q\in\Q$. 
So we  need to show that  
$$\alpha(\gamma_0)=\lambda_{i_0}(\gamma_0)=\pi (2-\lambda_{i_0}(\tau h_{\lambda_{i_0}}))=\pi (2-\lambda_{i_0}( 
h_{\tau\lambda_{i_0}}))\not=0.$$
 Suppose by contradiction that   $\lambda_{i_0}( h_{\tau\lambda_{i_0}})=2$. Then
both $\lambda_{i_0}-\tau\lambda_{i_0}$ and $\lambda_{i_0}-2\tau\lambda_{i_0}$ are roots in $\Delta$.
Since by assumption $\lambda_{i_0}$ and $\tau\lambda_{i_0}$ have non-zero restrictions to $\a$,   it means that
$${1\over 2}(\lambda_{i_0}-\tau\lambda_{i_0})|\a=\alpha, \quad (\lambda_{i_0}-\tau\lambda_{i_0})|\a=2\alpha,\quad 
(\lambda_{i_0}-2\tau\lambda_{i_0})|\a\not=0, \alpha, 2\alpha,\eqno(5.1)$$
which is absurd. This concludes the proof of (ii).
\sn
(iii)
 Assume now  that $\Delta_\a$ is reduced. Formulas (5.1) imply that $\lambda_{i_0}( h_{\tau\lambda_{i_0}})\le 0$. If 
$\lambda_{i_0}(\tau 
h_{\lambda_{i_0}})<0$, one has that $\tau\lambda= -\lambda$ (see [OS], Lemma 3.10). Then
$\lambda_{i_0}(h_{\tau \lambda_{i_0}})=-2$ and $\lambda_{i_0}(\gamma_0)=4\pi$.
\pn
If $\lambda_{i_0}(\tau h_{\lambda_{i_0}})=0$,  one has that $\lambda_{i_0}(\gamma_0)=2\pi$. To conclude the proof, we need to show that
in both  cases, $\gamma_0$ is the shortest vector in $ \a\cap i{1\over 2}\Gamma$.
By the classification results in [OS], Sect. (5.8), the restricted root system $\Delta_\a$ is reduced if and only if  
$(\g,\h)=(\s\o(p+1,q+1),\s\o(p+1,q ) $.  In this case, the root system $\Delta$ is either of type $B_n$ or of type $D_n$, depending on 
whether $ p+q$ is odd or even.
In the first case, any root $\lambda_{i_0}$ whose restriction to $\a$ is equal to $\alpha$ is a long root ([OS], p.466, [Wa], p.30). In 
the second case, all roots in $\Delta$ have the same length. By Lemma 5.1 and Lemma 5.2, in both cases $\gamma_0$ is the shortest 
vector 
in $ \a\cap i{1\over 2}\Gamma$ and the proof of (iii) is complete.

\bn
\proclaim Proposition 5.4. Let $G/H$ be a rank-one symmetric space. 
Then the polar map $\phi\colon \Omega \longrightarrow G^\C/H^\C$ is globally injective.

\mn\pf
Let $[g_1,X],[g_2,Y]\in\Omega$ be two points with the same image  $\phi([g_1,X])=\phi([g_2,Y])$.This means  that
$$g_1\exp iX=g_2\exp iY h_c,\quad\hbox{for some}~ h_c\in H^\C,~\hbox{and}~g_1,g_2\in G,$$ 
or equivalently 
$$  \exp iX= g \exp iY h_c \quad\hbox{for  }~ h_c\in H^\C,~\hbox{and}~g=g_1^{-1}g_2\in G.\eqno(5.2)$$
We want to show that $[g_1,X]=[g_2,Y]$, i.e. that there exists $h\in H$ such that 
$$\cases{ g_1=g_2h^{-1}\cr X=Ad_hY.}$$
Let $X=X_s+X_n$ and  $Y=Y_s+Y_n$  be the Jordan decomposition of $X$ and $Y$ in $\q$. Write 
$x=\exp iX= x_sx_n$, with $x_s=\exp iX_s$ and $x_n=\exp i X_n$, and similarly  
$y=\exp iY= y_sy_n$, with $y_s=\exp iY_s$ and $y_n=\exp i Y_n$.  
By Remark 4.1  and [Ma], Prop.2 (ii)a, p.66, 
equation (5.2) is equivalent to the system
$$     \cases{\exp iX_n= g\exp iY_ng^{-1}\cr \exp iX_s=g\exp iY_sh_c.} \eqno(5.3)$$
Let $\omega\subset \q$ be the subset defined in (2.7).
Observe that in the rank-one case, for every  non-zero semisimple element $Z\in \omega$, the element $z=\exp iZ\in\G$ is necessarily 
regular semisimple with respect to $\sigma, \tau$. In particular $Z$  is regular semisimple in $\q$. 
As a consequence,  the elements $X,Y$ in (5.2) are either both nilpotent or both semisimple.  
 In the first case,  system (5.3) reduces to
$$  \cases{\exp iX_n= g\exp iY_ng^{-1}\cr  g h_c=e.}  $$
Since the exponential map is injective on the set of
nilpotent elements, from the equation $\exp iX_n= g\exp iY_ng^{-1}=  \exp iAd_g  Y_n $,  we obtain 
$$X_n=Ad_g Y_n,\quad \hbox{with}~g\in  G\cap H^\C=H,$$ 
as requested.
If $X,Y$ are both  semisimple, system (5.3) reduces to the equation 
$$\exp iX= g\exp iY h_c \quad g\in G,~h_c\in H^\C. \eqno(5.4)$$    
Moreover, by [Ma], Thm.3,  the elements $X,Y $ may be assumed to sit in the same Cartan subspace in $ \q$. 
By Lemma 4.3(iii), equation (5.4)  implies then
$$\exp i4X= g\exp i4Yg^{-1}=\exp i4Ad_gY.\eqno(5.5)$$ 
Recall that  the compact Cartan subspace $\t$ of $\q$ is all contained in $\omega$ and that the restriction of the exponential map 
$\exp\colon  i\t\rightarrow \G$ is injective.
Hence, if $X,Y\in\t$, equation (5.5) implies 
$$X=Ad_gY,\quad \hbox{and}\quad x= \exp iX=g\exp iY g^{-1}=gyg^{-1}.$$ From the  last relation  one obtains that
 $g^{-1}=h_c\in G\cap 
H^\C=H$, as requested.
\sn
Assume now that $X,Y\in\omega\cap \a$, where $\a$ is the non-compact Cartan subspace in $\q$.  
By Lemma 4.3(v), one has that  $Ad_{g^{-1}}X\in\a$ and $g\in N_G(\a)$. 
Moreover, by Lemma 4.3(iv)  one has that $g^{-1}\tau(g)\in Z_H(\a)$. Therefore from (5.2) and Lemma 4.3(ii) it follows that
$$y^2 =g^{-1}x^2\tau(g)=g^{-1}x^2g c,\quad \hbox{for}~ c=g^{-1}\tau(g) \in Z_H(\a),~c^2=1. $$
Write the element $c=g^{-1}x^{-2}gy^2\in \exp i\a\cap \{x^2=1\}$ as
$$ c=\exp i \gamma ,\quad\hbox{with}~\gamma= 2(Y-Ad_{g^{-1}}X)\in\a\cap {1\over 2}i\Gamma.$$
Let $\alpha$ be the simple root $\alpha\in\Delta_\a$. Since $X,Y\in \a\cap\omega$, it follows from (2.7) that   
$$\cases{ |\alpha( \gamma) |=|\alpha( 2(Y-Ad_gX))|<2\pi, \hbox{in the reduced case}\cr 
 |\alpha( \gamma) |=|\alpha( 2(Y-Ad_gX))|<\pi, \hbox{in the non-reduced case.}\cr}$$
By  Lemma 5.3,  the element $c$ is  necessarily the identity element in $G$ and
$$y^2= gx^2g^{-1}.\eqno(5.6)$$
By  Lemma 4.3(ii), condition (5.6)  implies that $g\in H$. Finally, recall that $2\omega$ is contained in the injectivity set of
$\exp\colon \g^\C\rightarrow \G$. Then from (5.6)  or equivalently from 
$\exp i2Y=\exp i2Ad_gX$, it follows that
$Y=Ad_gX,~g\in H$, as requested.

\bn
The preceeding proposition shows that in the rank-one case, the domain $D=\phi(\Omega)$ is an analogue of the 
Akhiezer-Gindikin domain. The push-forward of the canonical pseudo-K\"ahler metric  of $\Omega$, deriving from the adapted complex 
structure, defines a 
$G$-invariant  pseudo-K\"ahler metric on $D$, of the same signature $(\sigma^+,\sigma^-)$ as the metric on $G/H$.
Next  we   show that $D$ is an increasing union   
  of $q$-complete smoothly bounded domains, where $q= \sigma^-$. As we shall see in the examples at the end of this section, in general 
$D$ is  not a Stein domain. 
\pn
(Recall that an $n$-dimensional  complex manifold $X$ is called $q$-complete if it admits an exhaustion function of class $C^2$,  whose 
 complex Hessian has 
at least $n-q$ non-negative  eigenvalues. By this definition, a Stein manifold is $0$-complete).

\bn 
Let $E\colon \Omega\longrightarrow \R$ be the Energy function defined in (3.1). Since $E$ is $G$-invariant, one has that $E([g,X])={1\over 2}B_\g(X ,X )$, where  $B_\g$ denotes the restriction of the Killing form 
of $\g$ to $\q$. It follows that  $E([g,N])=0$, for every nilpotent element $N\in\q$. Let $\t$ be the compact Cartan subspace in $\q$ and let $\Delta_\t=\Delta_\t(\g^\C,\t^\C)$ denote the corresponding restricted root system. Then $\t\subset \omega$  and for $X\in\t $, one has 
$$ E([g,X])= -(\RIm\alpha(X))^2(\dim \g^\alpha+4\dim \g^{2\alpha})\le 0,\quad \alpha\in \Delta_\t, ~   (\g^{2\alpha}~\hbox{possibly trivial}).$$
In particular, $E$ is non-positive and   $ E([g,X])\rightarrow  -\infty$,  for $|\RIm\alpha(X)|\rightarrow  \infty$. Let $\a$ be the non-compact Cartan subspace in $\q$ and let $\Delta_\a=\Delta_\a(\g^\C,\a^\C)$ denote the corresponding restricted root system.
For  $X\in\omega\cap\a$, one 
has 
$$ E([g,X])= \alpha(X)^2(\dim \g^\alpha+4\dim \g^{2\alpha}) ,\quad \alpha\in \Delta_\a, ~   (\g^{2\alpha}~\hbox{possibly trivial}).$$
Let $S$ be the supremum of $E$ on $\Omega$. 
For $0<s<S$,  define
$$\Omega_s=\{[g,X]\in \Omega~|~ E([g,X])<s\},\quad\hbox{and}  \quad D_s=\phi(\Omega_s).$$

\bn
\proclaim Proposition 5.5.   For every $s\in~]0,S[$, the $G$-invariant domain  $ D_s $ is   
 $q$-complete, for $q= \sigma^- $.   The domain $D$ is an increasing union of   $q$-complete domains
$$D=\bigcup_{0<s<S} D_s.$$

\mn\pf It is clear that $D=\bigcup_{0<s<S} D_s.$ It remains to show that each domain $D_s$ is $q$-complete, for $q=\sigma^-$.
For every $s\in~]0,S[$, the  boundary
$\partial D_s$ is a  regular orientable hypersurface. It consists of one or two  closed hypersurface $G$-orbits  intersecting  the  
slice 
$A=\exp i\a$. 
We  can compute the  Levi form  of the boundary of $D_s$ by computing the Levi form of these orbits.  Let $x_0=\exp iX_0\in\partial D_s$ be a base point, with $X_0\in\a\cap\omega$.  By [Ge], Prop. 5.14(i),  the Levi form of $\partial D_s$ at $x_0$ is a Hermitian matrix whose coefficients, up to a positive scalar multiple, are given by
$$L(\partial  D_s)_{x_0}\sim\pmatrix{I_{m^+(\alpha)}&0&0&0\cr 0&-I_{m^-(\alpha)}&0&0\cr 0&0&I_{m^+(2\alpha)}&0\cr 
0&0&0&-I_{m^-(2\alpha)}}.$$ 
Here the numbers $ m^+(\alpha),~m^-(\alpha),~m^+(2\alpha),~m^-(2\alpha) $  are the dimensions of the $\pm 1$-eigenspaces of the 
involution $\tau\theta$ on the root spaces $\g^\alpha$ and $\g^{2\alpha}$.
They are called the ``signatures" of the restricted root spaces (cf. [OS]).  
In our case, the numbers $ m^+(\alpha),~m^-(\alpha),~m^+(2\alpha),~m^-(2\alpha) $ are  given by
 $$ \m_+(\alpha)+m_+(2\alpha)=\dim\q\cap\p -1= \sigma^+-1,\qquad 
m_-(\alpha)+m_-(2\alpha) =\dim\q\cap\k= \sigma^-.$$  
Observe that the Levi form is positive definite when $\tau=\theta$ and $\dim\q\cap\k=0$. 
By [EVS], Thm. 3.8, p.421, a smoothly bounded open set in a Stein manifold, satisfying the above conditions,  is $q$-complete, for 
$q=\sigma^-$.  

\bn
{\bf Remark 5.6.}  The boundary $\partial D$ of $D$ is not  smooth. In the examples below, the Levi form of $\partial D$ at the smooth 
points is indefinite and degenerate. This shows that in general $D$ is not a Stein domain. In all such examples the manifold $\G/H^\C$ 
contains  no  $G$-invariant Stein subdomains. 

\bn
{\bf Example 5.7.}  {\it The real hyperboloids}. 
\pn
Let $p,q $ be positive integers, $p,q>2$. In $\C^{p+q}$ consider the manifold 
$$\X=\{Z \in \C^{p+q}~|~ z_1^2+\ldots+z_p^2-z_{p+1}^2-\ldots-z_{p+q}^2=-1\},\qquad \dim_\C\X=p+q-1 .$$
The group $\G=SO(p,q,\C)$ acts  transitively on $\X$. Taking as a base point 
$\xx=(0,\ldots,0,1)$,  there is an identification
$$\X=G^\C/H^\C=SO(p,q,\C)/SO(p,q-1,\C) .$$ 
Consider on $\X$ the action  of the connected real form $G=SO(p,q)_0$.
 It turns out that there are two pseudo-Riemannian $G$-symmetric spaces, embedded in $\X$ as totally real submanifolds of maximal 
dimension. To each of them there is associated a domain, image of the corresponding polar map. We determine such domains and examine 
their complex analytic properties.
 \sn
We begin by describing the $G$-orbit structure of $\X$. The  $G$-orbits in $\X$ are in one-to-one correspondence with the following set
\bn
$$\matrix{& & &\nn & &  \mm& & &&\cr 
& &Q(s)& *&P(t)& *& R(\sigma)& &&\cr 
 &---- &\overline{ \quad \quad \quad \quad \quad \quad  }& \bullet&  \overline{ \quad \quad \quad \quad \quad \quad  }&\bullet & 
\overline{ \quad \quad \quad \quad \quad \quad  }&----\cr
 & & &{ G/H}& &G/L &  & & &\cr   } \eqno(5.7)$$
\bn
The left black dot corresponds to the $G$-orbit of the point 
 $\xx=(0,\ldots,0,1)$, diffeomorphic to the pseudo-Riemannian rank-one symmetric space  $G/H= SO(p,q )_0/SO(p,q-1)$, of signature 
$ (p,q-1)$;  
 the right black dot  corresponds to 
the $G$-orbit of the point $\yy=(i,0,\ldots,0)$, diffeomorphic to the pseudo-Riemannian rank-one symmetric space  
$G/L= SO(p,q )_0/SO(p-1,q )$, of signature $ (p-1,q)$. Both $G/H$ 
and $G/L$ are  totally real submanifolds of $\X$, of real dimension 
 $$\dim_\R G/H=\dim_\R G/L=\dim_\C \X.$$
The central segment and the two halflines parametrize the  slices meeting the three types of closed orbits of maximal dimension. 
Since $G/H$ and $G/L$  has rank one, closed orbits of maximal dimension are real hypersurfaces in $\X$.
Let $\g=\h\oplus \q$ be the symmetric algebra  associated to $G/H$ and  $\g=\l\oplus \m$  the symmetric algebra  associated to  $G/L$.   
The  points $$Q(s)=(0,\ldots,0,i\sinh s,\cosh s), \quad \hbox{for}~s\in~]-\infty,0[,$$
parametrize the orbits intersecting the slice $C=\exp i\t \cdot \xx$, where $\t $ is the compact Cartan subspace in $\q$;
the  points
$$P(t)=(i\sin t ,0,\ldots,0,\cos  t),\quad  \hbox{for}~t\in ]0,\pi/2[ ,$$ parametrize the orbits intersecting the slice $A=\exp i\a $, 
where $\a $  is 
the non-compact Cartan subspace both in $ \q$ and in $\m$; 
the  points  
$$R(\sigma)=(i\cosh \sigma,\sinh \sigma, 0,\ldots,0),\quad \hbox{for}~\sigma\in~]0,\infty[,$$  parametrize the orbits intersecting the 
slice $C'=\exp i\t' \cdot \yy$, where $\t' $ is the compact Cartan subspace in $\m$.
Observe that 
$$Q(0)=P(0)=\xx,\qquad P({\pi\over 2})=Q(0)=\yy.$$
In addition to the closed orbits, there are two nilpotent hypersurface orbits: one with base point  
$$\nn=(i,0,\ldots,0_p,i,0,\ldots,1)$$ having $G/H$ in its closure, and the other with base  point
$$\mm=(i,0,\ldots,1_p,1,0,\ldots,0),$$
  having $G/L$ in its closure (here the 
subscript $p$ marks the $p^{th}$ coordinate of a point).
\bn
The  domain $D$ in $G^\C/H^\C$ determined by $G/H$ is given by 
$$D=~~\bigcup_{s<0}G\cdot Q(s)~~ \bigcup~~ G/H\cup G\cdot \nn~~ 
 \bigcup_{t\in ]0,{\pi\over 2}[} G\cdot  P(t),\qquad   \partial D=G/L \cup G\cdot \mm,$$
and corresponds to the following subset of diagram~(5.7)
$$\matrix{& &  &\nn& & \cr
& &Q(s)&* &P(t)& \cr 
 &---- &\overline{ \quad \quad \quad \quad \quad \quad  }&\bullet&\overline{ \quad \quad \quad \quad \quad \quad  }  \cr
 & & &{ G/H} &   &\cr   } $$
 \sn 
Consider now the real-valued $G$-invariant function
$$F\colon\X\longrightarrow \R,\quad Z\mapsto F(Z):= |z_1|^2+\ldots+|z_p|^2-|z_{p+1}|^2-\ldots-|z_{p+q}|^2.$$
Evaluating $F$ on the base points of all $G$-orbits in $\X$, it is easy to see that   $F$  separates the closed $G$-orbits  
$$-\infty\le F(Q_s)=-(\sinh^2 s+\cosh^2 s)< F(\xx)=-1< F(P_t)=1-2\cos^2 t< $$
$$< F(\yy)=1< F(R_\sigma)=\sinh^2 \sigma +\cosh^2 \sigma\le +\infty. \eqno(5.8)$$
The domain $D$ and its boundary $\partial D$ can be easily described by means of $F$:
$$D=\{Z\in\X~|~F(Z)<1\},\qquad \partial D=\{ Z\in\X~|~F(Z)=1\}.$$
The boundary  $\partial D$ is not smooth (the set of smooth points  coincides with the orbit $G\cdot \mm$).  The Levi form of $\partial D$ at the smooth points is degenerate  with $p-2$ positive 
eigenvalues, $q-1$ negative eigenvalues and 1 zero eigenvalue.  Hence, for $p,q$ sufficiently large, $D$ is not a Stein domain.
Let $R\in]-1,1[$.  For every such $R$, the subdomain  $D_R=\{Z\in D~|~F(Z)<R\}$  has smooth boundary,  with non-degenerate  Levi form of 
signature 
$(p-1,q-1)$. By  [EVS], Thm. 3.8, the domain $D_R$ is $(q-1)$-complete.
As a result, the domain $D$ is an increasing union of  $(q-1)$-complete domains 
$$D=\bigcup_{R<1}D_R.$$ 
\bn 
The  domain $D'$ in $G^\C/H^\C$ determined by $G/L$ is given by 
$$D'=~~\bigcup_{t\in ]0,{\pi\over 2}[}G\cdot P(t)~~ \bigcup~~ G/L\cup G\cdot \mm~~ 
 \bigcup_{\sigma>0} G\cdot  R(\sigma),\qquad   \partial D=G/H \cup G\cdot \nn,$$
and corresponds to the following subset of diagram~(5.7)
$$\matrix{  &  &\mm& && \cr
  &P(t) )&* &R(\sigma)&& \cr 
   &\overline{ \quad \quad \quad \quad \quad \quad  }&\bullet&\overline{ \quad \quad \quad \quad \quad \quad  }&----  \cr
   & &{ G/L} &   &&\cr   } $$
In terms of $F$,  the domain $D'$ and its boundary are given by
$$D'=\{Z\in\X~|~-F(Z)< 1\},\qquad \partial D'=\{ Z\in\X~|~F(Z)=-1\}.$$
The set of smooth points in $\partial D'$ coincides with the orbit $G\cdot \nn$ and has degenerate Levi form with $p-1$ positive 
eigenvalues, $q-2$ negative eigenvalues and 1 zero eigenvalue.  Hence, for $p,q$ sufficiently large, the domain $D'$ is not Stein 
either.
Let $R\in ]-1,1[$. For every such $R$,  the subdomain  $D'_R=\{Z\in D~|~-F(Z)<-R\}$ has smooth boundary  with non-degenerate Levi form of signature 
$(q-1,p-1)$. 
By [EVS], the domain $D'$ an the increasing union of $(p-1)$-complete domains: $D'=\cup_R D'_R$.
\bn
A few more remarks: for every $s\in ]-\infty,-1[$, the  level hypersurface 
$\X_R=\{F(Z)=R\}$ consists of one orbit with base point $Q_s=(0,\ldots,0,i\sinh s,\cosh s)$, for  $s\in~]-\infty,0[ $;  its Levi form   
is  non-degenerate with signature $(p,q-2)$. 
For every $R\in ]-1,+\infty[$, the  level hypersurface 
$\X_R=\{F(Z)=R\}$ consists of one orbit with base point $R_\sigma=(i\cosh \sigma,\sinh \sigma, 0,\ldots,0)$, for 
$\sigma\in~]0,+\infty[$; 
 its Levi form  is non-degenerate with  signature $(p-2,q)$. 
From these computations we  conclude that for $p,q$ sufficiently large,
the complex Hessian of every $G$-invariant function on $\X$ has both positive and negative eigenvalues at all points.  As a 
consequence, no $G$-invariant open set in $\X$ admits $G$-invariant plurisubharmonic functions.  
One can easily check that the $G$-action fails to be proper on every $G$-invariant open subset of $\X$. So no $G$-invariant open subset of $\X$ carries a $G$-invariant K\"ahler structure.

\beginsection 6. The higher rank case.

Let $G/H$ be a pseudo-Riemannian symmetric space  of rank higher than one. Let $\Omega\subset G\times_H\omega$ be the domain defined in 
(2.7). In this case, the polar map 
$$\phi\colon \Omega \longrightarrow  G^\C/H^\C$$
 is generally non-injective. As a result the domain $\Omega $  with the adapted complex structure  
 is a non-injective Riemann domain over $G^\C/H^\C$.

\sn
One may see this  as follows: 
let  $\c=\c_k\oplus \c_\p$ be a maximally split $\theta$-stable  Cartan subspace in $\q$, with both $\c_\k$ and $\c_\p$ different from 
$\{0\}$, and let $\Delta_\c=\Delta_\c(\g^\C,\c^\C)$ denote the corresponding restricted root system. Fix an element $X \in \c_\k$ with 
the 
property that
$\RIm \alpha(X )\not=0$, for all imaginary and complex roots in $\Delta_\c$. This is possible by taking $X$ in the complement of 
the finite set of hyperplanes $\{ H\in \c_\k~|~\RIm \alpha(H)=0\}_{\alpha\in\Delta_\c}$ in $\c_\k$. Let  $\gamma$ be an 
element in the intersection  $\c_\p\cap i\Gamma $,   where $\Gamma $ denotes the unit lattice in $\uu\subset \g^\C$. One such element  
$\gamma$ can be constructed as follows: 
let $\b\oplus \c$ be  a $\tau,\theta$-stable Cartan subalgebra in $\g$ extending $\c$. Let  $\lambda$ be a root in 
$\Delta(\g^\C,\b^\C\oplus \c^\C)$ with non-zero restriction to $\c$ and let $h_\lambda\in  \b_\R\oplus\c_\R$ be its inverse root (with 
$\lambda(h_\lambda)=2$). Then $\gamma=2\pi (h_\lambda+\theta h_\lambda+\tau h_\lambda+\theta\tau h_\lambda)$ lies  in
$\c_\p\cap i\Gamma $.
The elements $X$ and $X+\gamma$ satisfy the conditions 
$$X ,~X +\gamma\in \omega \quad {\rm and}\quad \exp i X =\exp i(X +\gamma).$$
As a consequence, the corresponding elements  $[e,X ]$ and $[e,X +\gamma]$ in $\Omega$ have the same image under the polar map.
Moreover, since  
the inclusion $ Z_H(X)\subset Z_H(x^2)$ may be a proper one, by Lemma 4.2(iv), the polar map $\phi$ may also fail to be 
 injective on some $G$-orbits in $\Omega$.

\bn
In the next proposition, we  show  that the polar map is injective on every closed orbit  of maximal dimension and is a covering map,
 when restricted to certain distinguished $G$-invariant subsets of $\Omega$. Such sets are coverings of  principal orbit strata in $D$. 
 \sn
Recall that when $G$ acts on $\G/H^\C$, closed orbits of maximal dimension come in a finite number of orbit types. Such orbits are 
called  
principal orbits, since  their union  is  an open dense subset of $\G/H^\C$.  
The set of principal orbits of a given type  is an open subset in  $\G/H^\C$, generally disconnected. It is referred to 
as a principal orbit stratum. The domain $D=\phi(\Omega)\subset \G/H^\C$ contains a number of connected components of principal orbit 
strata.
To be more precise,  
let $\c$ be a Cartan subspace of $\q$. 
Denote by $\c_{rs}$ the set of all elements $X\in\c$ with the property that $x=\exp iX\in\G$ is a regular semisimple element with 
respect 
to $\sigma,\tau$ (cf.(4.1)). The set $G  \exp i\c_{rs}H^\C$ is an open subset of $\G/H^\C$, consisting of closed orbits of maximal 
dimension, all of the same type. One of its connected components is contained in  $D=\phi(\Omega)$. 

\sn
Consider now the subset $$\Omega_\c :=G\times_{N_H(\c)}\c_{rs}\quad \subset\quad G\times_H\q.$$ 
Observe that  $\c_{rs}$ is  stable under the group $N_H(\c)$, so the set $\Omega_\c$ is well defined and the restriction of the 
polar map $\phi$ to $\Omega_\c$ has non-singular differential. 

\proclaim Proposition 6.1. The restriction of the polar map to $\Omega_\c$
$$\phi\colon \Omega_\c \longrightarrow \G/H^\C$$ is a $G$-equivariant covering map.

\mn\pf  
First we   prove that the restriction of
$\phi$ to every $G$-orbit in  $ \Omega_\c$ is injective.  Let  $[e,X]\in \Omega_\c $ and let $\bar x=\phi([e,X])$ be its image in 
$\G/H^\C$. Since the map $\phi$ is equivariant, it induces an inclusion of the corresponding isotropy subgroups 
$G_{[e,X]}\hookrightarrow 
G_{\bar
x}$.  The injectivity of $\phi$
on the orbit $G\cdot [e,X]$ is equivalent to showing that $G_{[e,X]}=G_{\bar x}=Z_H(X)$. This  follows  from Lemma 4.2(v).
Since the map
$$  \c_{rs}\longrightarrow \exp i\c_{rs} \longrightarrow \exp i\c_{rs}/ \exp i\c_{rs}\cap H^\C$$
is a covering map,  the proof of the proposition is complete.

\bn
In the next proposition, we show that the polar map is always injective on a  smaller $G$-invariant open subdomain $\Omega'\subset 
\Omega$, defined by 
$$ \Omega' =G\times_H \omega' ,\qquad \omega'=\{X\in\q~|~|\RRe\lambda|<\pi/4,~\hbox{ for all}~ \lambda\in spec(ad_X)\}.\eqno(6.1)$$

\proclaim Proposition 6.2.  The polar map $\phi$ is  injective on the domain $\Omega'=G\times_H \omega'  $.

\mn\pf Let $[g_1, X ],[g_2, Y]\in G \times_H\omega' $ be points with the same image in $\G/H^\C$, i.e.
$$ g \exp iX =\exp iY h ,\quad\hbox{for}~  h \in H^\C,~g=g_2^{-1}g_1\in G.\eqno(6.2)$$  
Lemma 4.3(iii) implies that
$$\exp i 4Y=g \exp i 4X g^{-1}=\exp i Ad_g 4X .$$
Moreover, since $4\omega'$ is contained in the injectivity set of $\exp\colon \g^\C\rightarrow \G$, one has that
$Y=Ad_g X $. This relation together with  equation (6.2) implies that
 $g=h\in G\cap H^\C=H$. In other words
$$\cases{g_2=g_1h^{-1}\cr Y=Ad_{h}X },\quad h\in H,$$
as requested.

\bn
{\bf Corollary 6.3.} By Proposition 6.1, a canonical $G$-invariant pseudo-K\"ahler metric is defined on coverings of principal orbit 
strata in $D$. 
Proposition 6.1 and 6.2 extend similar results obtained in [Fe] and [Br2] for the group case, by  different methods.

\bn
{\bf Remark 6.4.}  {\it If $G/H$ is a non-Riemannian semisimple symmetric space,  the domains  $D=\phi(\Omega)$ and  $D'=\phi(\Omega')$ 
cannot be hyperbolic}.
Both  $D $  and $D'  $ contain  the complex homogeneous subvariety 
$$K^\C/(K\cap H)^\C\quad \hookrightarrow \quad D' ~~\subset ~~D, $$
embedded as the $K^\C$-orbit of the base point $eH^\C$.
The homogeneous subvariety $K^\C/(K\cap H)^\C$ is also the image  of the set $K\times_{K\cap H}\k\cap\q\subset \Omega'\subset \Omega$ 
by the map $\phi$. Indeed, for every $X\in \k\cap \q$, the eigenvalues of $ad_X$ are all purely imaginary and the restriction of   
$\phi$ to $K\times_{K\cap H}\k\cap \q$ has both non-singular differential and it is injective.
It can be viewed as the complexification of the compact symmetric space $K/K\cap H\hookrightarrow G/H$, embedded in $G/H$ as the 
$K$-orbit of the base 
point $eH$. One has that
$$ \dim_\C  K^\C/(K\cap H)^\C=\dim_\R K/K\cap H=\dim_\R \q\cap \k.$$

\bn\bn\bn
\centerline{{\bf References}}
\bn\bn 
\item{[AG]}  Akhiezer D., Gindikin S., {\sl On Stein extensions of real symmetric spaces}, Math. Ann., 
{\bf 286} (1990), 1--12. 
\smallskip 
\item{[Ba]} Barchini L., {\sl Stein extensions of real symmetric spaces and the geometry of the flag manifold}, 
Math. Ann., {\bf 326} (2003) 331--346.
\smallskip 
\item{[BHH]} Burns D., Halverscheid S., Hind R., {\sl The Geometry of Grauert Tubes and Complexification of Symmetric Spaces}, Duke 
Math. J., {\bf 118} (2003) 465--491.

\smallskip 
\item{[Br1]} Bremigan R., {\sl Quotients for algebraic group actions over non-algebraically closed fields},
J. reine angew. Math. {\bf 453} (1994) 21--47. 
\smallskip 
\item{[Br2]}  Bremigan  R., {\sl Pseudo-K\"ahler forms on complex Lie groups}, Documenta Math. {\bf 5} (2000), 595--611.
\smallskip 
\item{[EVS]} Eastwood M., Vigna-Suria G., {\sl Cohomologically complete and pseudoconvex domains}, Comm. Math. Helv., {\bf 55} (1980) 
413--426.  
\smallskip 
\item{[Fe]}  Fels  G., {\sl Pseudo-K\"ahlerian structure on domains over a complex semisimple Lie group},
Math. Ann. {\bf 232} (2002), 1--29.
\smallskip 
\item{[H]} Huckleberry A.,  {\sl On certain domains in cycle spaces of flag manifolds}, 
Math. Ann., {\bf 323} (2002) 797--810.
\smallskip 
\item{[Ha]} Halverscheid S.,  {\sl On Maximal Domains of Definition of Adapted Complex Structures for Symmetric Spaces of Non-compact 
Type,} Thesis  Ruhr-Universit\"at Bochum, 2001.
\smallskip 
\item{[HI]} Halverscheid, S., Iannuzzi, A.,  {\sl Maximal complexifications of certain Riemannian
homogeneous spaces}. Trans. of the AMS,  {\bf 355}, 11 (2003) 4581--4594.

\smallskip 
\item{[HO]} Hilgert J., \`Olafsson G., Causal symmetric spaces. Geometry and Harmonic Analysis, Perspectives in Math. {\bf 18}, 
Academic 
Press, San Diego, CA, 1997.
\smallskip 
\item{[Hu]} Humphreys J.E., Conjugacy Classes in Semisimple Algebraic groups, AMS Math. Surveys and Monographs, Vol. 43, Providence RI, 
1995.
\smallskip 
\item{[KS1]} Kr\"otz  B., Stanton  R., {\sl Holomorphic extensions of representations I: Automorphic functions},  Ann. of Math. {\bf 
159}, 
2 (2004), 641--724.  
\smallskip 
\item{[KS2]} Kr\" otz B., Stanton R., {\sl Holomorphic extensions of representations II: geometry and harmonic analysis}, GAFA, Geom. 
Funct. Anal., {\bf 15} (2005), 190-245.
\smallskip 
\item{[Ge]} Geatti L., {\sl Invariant domains in the complexification of a non-compact
        Riemannian  symmetric space}, J. of Algebra, {\bf 251} (2002) 619--685.
\smallskip 
\item{[GS1]} Guillemin V., Stenzel M.,{\sl Grauert tubes and the homogeneous Monge-Amp\`ere equation}, J. Diff. Geom. {\bf 34} (1991) 
561--570.
\smallskip 
\item{[GS2]} Guillemin V., Stenzel M.,{\sl Grauert tubes and the homogeneous Monge-Amp\`ere
equation}, J. Diff. Geom. {\bf 35} (1992) 627--641.  
\smallskip
\item{[LS]} Lempert L., Sz\H oke R.,{\sl Global solutions of the homogeneous complex 
Monge-Amp\`ere equation and complex structures on the tangent bundles of Riemannian manifolds}, 
Math. Ann. {\bf 290} (1991) 689--712.
\smallskip
\item{[Ma]} Matsuki T., {\sl Double coset decomposition of reductive Lie groups 
arising from two involutions}, J. of Algebra, {\bf 197} (1997) 49--91.
\smallskip
\smallskip 
\item{[Ne1]} Neeb K.H., {\sl On the complex geometry of invariant domains in 
complexified symmetric spaces}, Ann. Inst. Fourier,
Grenoble, {\bf 49}, 1 (1999), 177--225. 
\smallskip
\item{[Ne2]} Neeb K.H., Holomorphy and convexity in Lie theory, de Gruyter Expositions in Mathematics  28, 
Walter de Gruyter  Co., Berlin, 2000.
\smallskip

\item{[OS]} Oshima T., Sekiguchi J., {\sl The restricted root system of a semisimple symmetric pair}, Advanced Studies in Pure Math. 4, 
Group Representations and Systems of Differential Equations, 1984, pp. 433--497.
\smallskip
\item{[St]}  Steinberg R., {\sl Torsion in reductive groups}, Adv. in Math., {\bf 15} (1975), 13--92.
\smallskip

\item{[Sz1]}  Sz\H oke R.,{\sl  Complex structures on the tangent bundles of Riemannian manifolds}, 
Math. Ann. {\bf 291} (1991) 409--428.
\smallskip

\item{[Sz2]}  Sz\H oke R.,{\sl Canonical complex structures associated to connections and complexifications of  Lie groups}, 
Math. Ann. {\bf 329} (2004) 553--591.
 \smallskip 
 \item{[vD]} van Dijk G., {\sl Orbits of real affine symmetric spaces I: the infinitesimal case},
Indag. Math. {\bf 45} (1983) 51--66. 
\smallskip 

\item{[Va]} Varadarajan V.S., Lie groups, Lie algebras and their representations, Springer Verlag, New York 1984.
\smallskip 
\item{[Wa]} Warner  G., Harmonic analysis on Semi-Simple Lie groups I, Springer Verlag, Berlin-Heidelberg-New York 1972.

\bye